\documentclass{article}

\usepackage{arxiv}

\usepackage[utf8]{inputenc} 
\usepackage[T1]{fontenc}    
\usepackage{hyperref}       
\usepackage{url}            
\usepackage{booktabs}       
\usepackage{amsfonts}       
\usepackage{nicefrac}       
\usepackage{microtype}      
\usepackage{lipsum}		    
\usepackage{graphicx}
\usepackage[square,sort,comma,numbers]{natbib}
\usepackage{doi}
\usepackage{amsmath}
\usepackage{subcaption}

\title{A spectral element solution of the 2D linearized potential flow radiation problem}

\date{April 29, 2022}	

\author{ \href{https://orcid.org/0000-0001-6698-2623}{\includegraphics[scale=0.06]{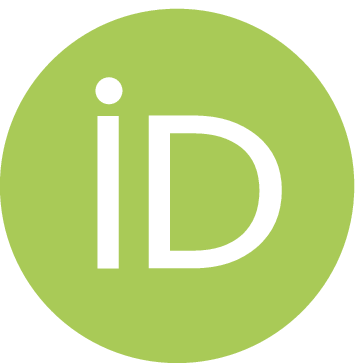}\hspace{1mm}Jens Visbech} \\
	Department of Applied Mathematics and Computer Science\\
	Technical University of Denmark\\
	Kongens Lyngby, 2800 \\
	\texttt{jenshaakon@gmail.com} \\
	\And
	\href{https://orcid.org/0000-0001-8626-1575}{\includegraphics[scale=0.06]{Figures/logo_orcid.eps}\hspace{1mm}Allan P. Engsig-Karup} \\
	Department of Applied Mathematics and Computer Science\\
	Technical University of Denmark\\
	Kongens Lyngby, 2800 \\
	\texttt{apek.dtu.dk} \\
	\And
	\href{https://orcid.org/0000-0002-7263-442X}{\includegraphics[scale=0.06]{Figures/logo_orcid.eps}\hspace{1mm}Harry B. Bingham} \\
	Department of Civil \& Mechanical Engineering\\
	Technical University of Denmark\\
	Kongens Lyngby, 2800 \\
	\texttt{hbb@mek.dtu.dk} \\
}

\renewcommand{\shorttitle}{Visbech et al. (2022)}

\hypersetup{
pdftitle={2D_Pseudo_Impulsive_SEM},
pdfsubject={Preprint},
pdfauthor={Jens Visbech, Allan P. Engsig-Karup, Harry, B. Bingham},
pdfkeywords={Spectral Element Method, Potential Flow, Water Waves, Pseudo-impulsive Formulation, High-order Numerical Method, Offshore Hydrodynamics},
}

\begin{document}
\maketitle

\begin{abstract}
 We present a scalable 2D Galerkin spectral element method solution to the linearized potential flow radiation problem for wave induced forcing of a floating offshore structure. The pseudo-impulsive formulation of the problem is solved in the time-domain using a Gaussian displacement signal tailored to the discrete resolution. The added mass and damping coefficients are then obtained via Fourier transformation. The spectral element method is used to discretize the spatial fluid domain, whereas the classical explicit 4-stage 4th order Runge-Kutta scheme is employed for the temporal integration. Spectral convergence of the proposed model is established for both affine and curvilinear elements, and the computational effort is shown to scale with $\mathcal{O}(N^p)$, with $N$ begin the total number of grid points and $p \approx 1$. Temporal stability properties, caused by the spatial resolution, are considered to ensure a stable model. The solver is used to compute the hydrodynamic coefficients for several floating bodies and compare against known public benchmark results. The results are showing excellent agreement, ultimately validating the solver and emphasizing the geometrical flexibility and high accuracy and efficiency of the proposed solver strategy. Lastly, an extensive investigation of non-resolved energy from the pseudo-impulse is carried out to characterise the induced spurious oscillations of the free surface quantities leading to a verification of a proposal on how to efficiently and accurately calculate added mass and damping coefficients in pseudo-impulsive solvers. 
\end{abstract}


\section{Introduction}

The offshore industry is of considerable economic and ecological importance in the world. In particular, floating offshore structures play a vital role in keeping up with society's present and future demand for sustainable and green energy resources. Here will e.g. floating wind turbine foundations allow for the establishment of energy farms in greater water depths, where bottom-mounted foundations are too expensive to install \cite{andersen2016floating}. Further development of wave energy converters, floating solar panel plants, and seaweed farms can also contribute in reducing the need for fossil-based fuels for energy production purposes \cite{esmaeilzadeh2019shape}.

\subsection{Choosing the mathematical formulation}

Designing structures at sea requires a very detailed, in-depth analysis of the physical environment to ensure a safe, well-functioning, and sustainable result. The study of waves and the interaction between waves and floating structures can be done in different ways, as it has been for many decades. One way is to use numerical methods and mathematical formulations that mimic the physical environment. The most precise mathematical formulation - yet computationally expensive - is the well-known Navier-Stokes Equations (NSE). The complexity and computational expense can be reduced dramatically by applying various assumptions to these equations, though, at a cost in terms of physical accuracy. A recent comprehensive review of wave energy converter design and different mathematical formulations is given in Davidson and Costello (2020) \cite{davidson2020efficient}.

The nonlinear and linear Potential Flow (PF) simplification is widely used for free surface gravity waves in the field of offshore and marine engineering. Compared to the NSE, the PF formulations will not be able to simulate breaking waves, overturning waves, nor viscous effects e.g. at the seabed or near structures; however, in offshore environments, it is often sufficiently accurate to model these effects through empirical boundary terms. The PF formulations require only a fraction of the computational effort needed by NSE-based solvers, and are thus much more suitable for long-time and large-scale wave simulations \cite{ransley2019blind}. 

Considering wave-structure interaction and the linear PF setting, the classical impulsive time-domain formulation was widely used and studied in the 1980s and 90s \cite{korsmeyer1988first,king1988seakeeping,bingham1994simulating}. King (1987) \cite{king1987time} introduced the pseudo-impulsive formulation, and a significant amount of work on this topic has appeared over the past decade \cite{read2011linear,read2012overset,read2012solving,amini2017solving,amini2018pseudo}. The basic concept of impulsive formulations is to force the floating body with an impulse in velocity, let the system evolve in time, and measure the resultant force on the body. A detailed investigation of the difference between the classical impulsive and pseudo-impulsive formulations is given in Section \ref{sec:Investigation}.

\subsection{Choosing the numerical method for a potential flow formulation}

Wave-structure interaction problems have traditionally been solved numerically using the Boundary Element Method (BEM) where the formulation is recast as an integral equation, which is solved - in a discrete sense - on the domain boundaries. A key feature of the BEM is the ability to easily resolve complex geometries. With the projection onto the boundaries, the formulation is reduced by one dimension. This essentially means that every boundary node receives information from the entire fluid domain, which leads to a dense non-symmetric system of equations. With this making large-scale simulation challenging due to the $\mathcal{O}(N^2)$ computational work effort when solved using conventional direct or iterative solvers. By using multipole expansions and/or pre-corrected FFT methods, the system can be made much more sparse, enabling asymptotic $\mathcal{O}(N \log N)$ scaling, however with a large constant in front \cite{jiang2012precorrected,white1994comparing}. The BEM has been implemented into the widely used \textit{state of the art} frequency-domain-based commercial software WAMIT \cite{lee2006wamit}; and also in the open-source software NEMOH \cite{babarit2015theoretical}. Impulsive time-domain implementations have as well been made using the BEM, e.g. for zero speed problems \cite{korsmeyer1988first} and for forward speed problems \cite{bingham1994simulating,korsmeyer1998forward,kara2011prediction}. 

Approaches for solving the PF problem have also been developed using the Finite Difference Method (FDM), both for pure nonlinear wave propagation \cite{bingham2007accuracy,engsig2009efficient}, nonlinear wave-structure interaction \cite{ducrozet2010high}, linear wave-structure interaction using pseudo-impulsive time domain formulations at zero speed \cite{read2011linear,read2012solving,read2012overset}, and with forward speed \cite{amini2017solving,amini2018pseudo}. Using an FDM-based model, the entire fluid domain is discretized to solve the Laplace problem, hereby resulting in a larger linear system of equations compared to the BEM discretization; however, the system is sparse due to the locally applied stencils. This can - via iterative solvers - lead to an optimal $\mathcal{O}(N)$ scaling in work effort for higher order accurate schemes \cite{bingham2007accuracy,engsig2009efficient}. Complex geometries can also be represented using this discretization method; yet their inclusion is far from trivial since updating and adapting the mesh domain and imposing boundary conditions consistently can be challenging within this discretization framework.

The Finite Element Method (FEM) has also been used for modeling linear and nonlinear wave propagation and wave-structure interaction \cite{wu1994finite,taylor1996analysis,ma2001finite,wu2003coupled}. Due to the local support of the linear basis functions used in the FEM, the resulting linear system of equations will be sparse, hereby enabling optimal computational scaling and a small memory requirement. Also, the element discretization approach makes it very suitable for complex geometries; however, re-meshing strategies may be necessary in the nonlinear setting due to the movement of the free surface boundary. A drawback for the FEM is the restriction to second order in spatial accuracy from the use of linear basis functions.

This caveat naturally leads to another element based method, namely the higher order version - or extension - of the FEM, i.e. the spectral element method (SEM). This method can also be seen as a multi-domain version of the classical single domain spectral method \cite{canuto2007spectral,canuto2007spectral_evolution,kopriva2009implementing}, where higher order polynomial basis functions are used to enable higher order spatial accuracy. The use of SEM for fluid dynamics problems was originally proposed by Patera (1984) \cite{patera1984spectral}. In the late 90s, the SEM was shown to be inherently unstable for free surface flows using simplex elements by Robertson and Sherwin (1999) \cite{robertson1999free}; however, almost two decades later, the first stable SEM-based fully nonlinear PF model for pure wave propagation was successfully developed by Engsig-Karup et al. (2016) \cite{engsig2016stabilised}. Since this breakthrough, multiple extensions have been developed and applied to model wave-wave, wave-bed, and wave-structure interaction with optimal $\mathcal{O}(N)$ scaling through $p$-multigrid methods \cite{engsig2019mixed,laskowski2019modelling,engsig2021efficient,mortensen2021simulation}. A recent review of the \textit{state of the art} for SEM models can be found in Xu et al. (2018) \cite{xu2018spectral}. The clear advantage of the SEM is the geometrical flexibility, optimal computational scaling properties, and high spatial accuracy. In combination with the general increase in computational resources, long-time and large-scale simulations can be achieved using the SEM as the numerical discretization method for linear and nonlinear PF formulations.

\subsection{Paper contributions and overview}

In this work, we develop a linear potential flow spectral element (PF-SEM) solver for the simulation of linear wave-structure interaction through a pseudo-impulsive PF formulation. The unstructured mesh model has the ability to take into account the geometric complexity of the structures by using affine and curvilinear elements and is validated through convergence studies, stability considerations, and simulations that are compared against known benchmark results from the literature and \textit{state of the art} solvers. The computational effort is shown to scale with $\mathcal{O}(N^p)$, where $p \approx 1$. This confirms the potential of the model and completes and extends the preliminary results presented by Visbech et al. (2022) \cite{visbech2022efficient}. An additional novelty is the investigation of the non-physical behavior of the numerical solution at the free-surface/body intersection arising from unresolvable energy in the pseudo-impulse (or impulse), highlighting the importance of properly tuning the impulse signal.

The rest of the paper is organised as follows: In Section \ref{sec:Mathematical_Formulation}, we present the mathematical formulation in terms of the governing equations for the pseudo-impulsive formulation. Then, in Section \ref{sec:Numerical_Methods_and_Discretization}, we outline the numerical method and discretization, with this covering the SEM, temporal integration, computation of pressure and hydrodynamic coefficients, and the construction of the discrete pseudo-impulse. In Section \ref{sec:Numerical_Properties} the numerical properties are outlined by showing studies of convergence and computational scaling effort, but also, highlighting linear stability analysis of the spatial discretization. Section \ref{sec:Numerical_Experiments} seeks to display different numerical experiments obtained using the proposed model, and in Section \ref{sec:Conclusion}, a conclusion is made to summarize the novelty achievements and findings.
\section{Mathematical Formulation}\label{sec:Mathematical_Formulation}

We consider a spatial 2D $(d=2)$ fluid domain, $\Omega \subset {\rm I\!R} ^{d}$, upon which the flow is assumed to be: 1) impressible, 2) inviscid, and 3) irrotational. These simplifications enable the NSE to be represented via PF with the fluid velocity field, $\boldsymbol{u} = ( u(x,z,t),w(x,z,t) )$, expressed as the gradient of a scalar velocity potential, $\boldsymbol{u} = \nabla \phi$, where $\nabla$ is the spatial differential operator $\nabla = (\partial_x,\partial_z)$ in a Cartesian coordinate system spanned by the horizontal $x$-axis and the vertical $z$-axis. Let $\mathcal{T}$ denote the time domain by $\mathcal{T}:t \geq 0$, with $t$ being the time. Now $\phi = \phi(x,z,t):  \Omega \times \mathcal{T} \xrightarrow{} {\rm I\!R}$ and $\boldsymbol{u} = \boldsymbol{u}(x,z,t):  \Omega \times \mathcal{T} \xrightarrow{} {\rm I\!R}$.

An even further simplification of the PF equations is often used in offshore engineering by assuming that the waves and the oscillatory motions of the floating structure are of small amplitude compared to the wavelength and the size of the structure, respectively. This yields the linearized version of the PF equations, where the perturbation expansion term of first order of the complete nonlinear solution still captures the majority of the physics. The linear assumption enables the general solution for wave-structure interaction to be decomposed into independently solvable unsteady radiation and diffraction parts \cite{haskind1946hydrodynamic}, where the focus of this paper is on the former. Thus
\begin{equation}\label{eq:linear_decomposition_potential}
    \phi = \phi^R + \phi_0 + \phi_7,
\end{equation}
with $\phi_0$ denoting the incident wave potential and $\phi_7$ the scattered wave potential caused by the interaction of the fixed structure with the incident waves. Finally, $\phi^{R}$ denotes the total radiation potential, $\phi^{R} = \sum_{k = 1}^6 \phi_k$, due to rigid-body motions in six degrees of freedom, i.e. surge ($k=1$), sway ($k=2$), heave ($k=3$), roll ($k=4$), pitch ($k=5$), and yaw ($k=6$). For the 2D setting only surge, heave, and pitch are relevant.

\begin{figure}[t]
\begin{minipage}[c]{0.3\linewidth}
    \caption{Illustration of the 2D fluid domain, $\Omega$, and the associated boundaries, $\Gamma$. The free surface is denoted by $\eta(x,t)$ and the water depth by $h(x)$. The arbitrary domain length in $x$ is denoted by $L$.}
    \label{fig:2D_Domain}
\end{minipage}
\hfill
\begin{minipage}[c]{0.60\linewidth}
    \includegraphics{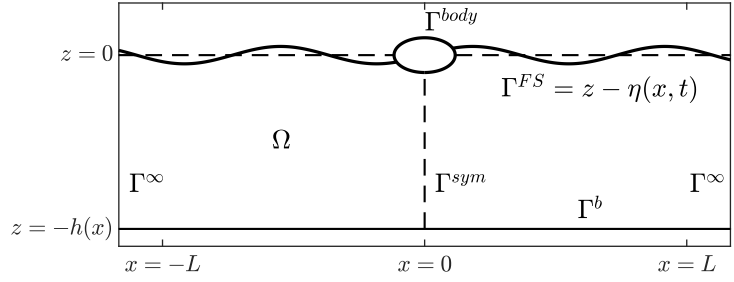}
\end{minipage}
\end{figure}

The entire fluid domain, $\Omega$, will be bounded by five different boundaries, $\Gamma = \left (\Gamma^{FS} \cup \Gamma^{b} \cup \Gamma^{\infty} \cup \Gamma^{body} \cup \Gamma^{sym} \right ) \subset {\rm I\!R}^{d-1}$, as illustrated on Figure \ref{fig:2D_Domain}, namely: 1) from above at $z=\eta \approx 0$ is the time-dependent free surface boundary, $\Gamma^{FS}$ , 2) from below at $z = -h$ by the bathymetry boundary, $\Gamma^{b}$, 3) from the ends of the domain by the far-field boundaries, $\Gamma^{\infty}$, 4) from the floating body by the body boundary, $\Gamma^{body}$, and in some special instances 5) from the symmetry boundary, $\Gamma^{sym}$. We have that $\eta = \eta(x,t):  \Gamma^{FS} \times \mathcal{T} \xrightarrow{} {\rm I\!R}$ is the free surface elevation and $h = h(x): \Gamma^{b} \xrightarrow{} {\rm I\!R}$ is the water depth. For the entirety of this work, the depth is assumed constant, $h(x) = h$.

\subsection{Governing equations}
The governing equations for each of the $k$ linear radiation potentials, $\phi_k$, is achieved using the pseudo-impulsive formulation by combining the linear Eulerian kinematic and dynamic free surface boundary conditions
\begin{align}\label{eq:free_surface_conditions}
    \partial_{t} \eta_k &= \partial_z \phi_k, \quad \text{and} \quad  \partial_{t} \phi_k = -g \eta_k, \quad \text{on} \quad \Gamma^{FS} \times \mathcal{T},
\end{align}
with the mass conserving Laplace boundary value problem including suitable boundary conditions as
\begin{align}\label{eq:Laplace_equation}
    \nabla^{2} \phi_k &= 0, \quad \text{in} \quad \Omega,\\ \label{eq:free_surface_BC}
    \phi_k &= \phi_k, \quad \text{on} \quad \Gamma^{FS},\\
    \partial_{n} \phi_k &=0, \quad \text{on} \quad \Gamma^b,\\
    \partial_{n} \phi_k &= \dot{x}_k n_k, \quad \text{on} \quad \Gamma^{body},\label{eq:body_BC}
\end{align}
where $g=9.81$ m/s$^2$ is the constant gravitational acceleration, $\partial_n$ is the derivative in the normal direction to the boundary, $\dot{x}_k = \dot{x}_k(t):\mathcal{T} \xrightarrow{} {\rm I\!R}$ is the pseudo-impulsive velocity from the Gaussian displacement, $x_k = x_k(t):\mathcal{T} \xrightarrow{} {\rm I\!R}$, that will be defined in Section \ref{sec:discrete_pseudo_impulse}, and $n_k$ is a generalized unit normal vector given as
\begin{equation}\label{eq:generalized_normal_vector}
n_{k}= \boldsymbol{n}, \quad \text{for} \quad k = \{1,3\}, \quad \text{and} \quad n_{k}=\boldsymbol{r} \times \boldsymbol{n}, \quad \text{for} \quad k = 5, \quad \text{on} \quad \Gamma^{body},
\end{equation}
where $\boldsymbol{n}$ is the unit normal vector and $\boldsymbol{r}$ is the position vector of a point on $\Gamma^{body}$.  

To complete the Laplace problem and make it well-posed, boundary conditions on the symmetry boundary, $\Gamma^{sym}$, and the far-field boundary, $\Gamma^{\infty}$, are needed. If $\Gamma^{sym}$ is applicable, it depends on the solution of the problem as
\begin{align}\label{eq:symmetri_BC}
\text{Anti-symmetric:} \quad \phi_k=0, \quad \text{and} \quad \text{Symmetric:} \quad  \partial_n \phi_k=0, \quad \text{on} \quad \Gamma^{sym},
\end{align}
where the former applies for $k=1$ and $k=5$, and the latter for $k=3$. On $\Gamma^{\infty}$ the main object is to ensure that no energy is reflected back into the computational domain, ultimately disturbing the solution. For this work, two different approaches are used to make a non-reflective domain: 1) An impermeability condition, $\partial_{n} \phi_k =0$, on $\Gamma^{\infty}$ combined with a horizontal grid stretching strategy \cite{xxx}, or 2) a classical relaxation zone \cite{larsen1983open} (for the high frequency content) combined with a Sommerfeld condition \cite{zhang2009non}, $\partial_{n} \phi_k = V_s$, on $\Gamma^{\infty}$ (for the low frequency content). Assuming a straight and vertical far-field boundary, the imposed velocity, $V_s$, at time $t=t_i$ is determined as $V_s = V_s(z,t_i) =  u(L-\Delta x,z, t_i - \Delta t)$, with $\Delta x = \Delta t \sqrt{gh}$ being determined from the shallow water limit. The discrete time step size, $\Delta t$, is defined in Section \ref{sec:Time_Integration}.

\subsection{The pseudo-impulsive radiation problem}

From each of the computed radiation velocity potentials, $\phi_{k}$, the dynamic pressure, $p_k = -\rho \partial_{t} \phi_{k}$, where $p_k=p_k(x,z,t): \Omega \times \mathcal{T} \xrightarrow{} {\rm I\!R}$, can be found from the linearized Bernoulli equation at every point in $\Omega$ for any time instance in $\mathcal{T}$. From this, the radiation forces, $F_{jk} = F_{jk}(t):  \mathcal{T} \xrightarrow{} {\rm I\!R}$, can be evaluated by integrating over the wetted-body surface, $\Gamma^{body}$. With this, the radiation force in the $j$'th direction due to a pseudo-impulsive motion in the $k$'th direction is
\begin{equation}\label{eq:F_jk}
    F_{jk}= \int_{\Gamma^{body}} p_k  n_{j} \ d \Gamma, \quad \text{for} \quad \mathcal{T},
\end{equation}
where $n_j$ is the generalized normal vector given by \eqref{eq:generalized_normal_vector}.

Let $\mathcal{H}$ denote the frequency space $\mathcal{H}:\omega \geq 0$, with $\omega$ being the radial frequency that is connected to the cyclic frequency as $f = \frac{\omega}{2 \pi}$. Following Amini-Afshar and Bingham (2017) \cite{amini2017solving}, the frequency-dependent added mass and damping coefficients, $(a_{jk},b_{jk}) = (a_{jk}(\omega),b_{jk}(\omega)): \mathcal{H} \xrightarrow{} {\rm I\!R}$, are given by
\begin{equation}\label{eq:a_jk_b_jk}
    \omega^2 a_{jk} - i \omega b_{jk} = \frac{\mathcal{F}\{F_{jk}\}}{\mathcal{F}\{x_k\}},
\end{equation}
where $\mathcal{F}\{\cdot\}$ denotes the Fourier transform operator, $\mathcal{F}:\mathcal{T}\xrightarrow{} \mathcal{H}$. 
\section{Numerical Methods and Discretization}\label{sec:Numerical_Methods_and_Discretization}

Exploiting a traditional Method Of Lines (MOL) approach, the governing equations are discretized spatially and kept continuous in time. We seek to represent the spatial part by the SEM, ultimately yielding a semi-discrete system of ordinary differential equations, which then are time-integrated using the classical Explicit 4-stage 4th order Runge-Kutta Method (ERK4). The time-derivative for the pressure computation needed in \eqref{eq:F_jk} will be performed by a 4th order finite difference approximation, and the Fourier transformations in \eqref{eq:a_jk_b_jk} are handled discretely using the Fast Fourier Transform (FFT)\cite{cooley1965algorithm}. Finally, the design of the discrete pseudo-impulsive velocity and displacement used in \eqref{eq:body_BC} and \eqref{eq:a_jk_b_jk}, respectively, will be specified. All of the aforementioned objectives are outlined in the following.

\subsection{Spatial discretization - The spectral element method}

The fluid domain, $\Omega$, will be tessellated using conforming non-overlapping affine or curvilinear triangular elements in an unstructured manner, where we can arrange the elements uniformly or apply mesh refinement at the free surface and body boundary, $\Gamma^{FS}$ and $\Gamma^{body}$. See Figures \ref{fig:First_three_meshes_affine}, \ref{fig:Mesh_curvilinear}, and \ref{fig:Mesh_Cylinder} for visualizations. The tessellation will consist of $N_{elm}$ elements such that we, at some particular time instance, have a time-invariant partitioning of the domain and the global radiation velocity potential as
\begin{equation}
    \Omega \simeq \bigcup_{n = 1}^{N_{elm}} E^n, \quad \text{and} \quad   \phi_k(x,z) = \bigoplus_{n = 1}^{N_{elm}} \phi_k^n(x^n,z^n), \quad \text{for} \quad \left(x^n,z^n \right) \in E^n, 
\end{equation}
with $E^n$ being the $n$'th element having $N_{ep}$ sets of coordinates $\left(x^n,z^n \right) = \left \{ \left(x^n,z^n \right) \right \}_{i=1}^{N_{ep}}$, and $\phi_k^n$ is the local solution defined on this element. We seek to apply the following approximate polynomial representation of the local solution
\begin{equation}\label{eq:polynomial_representation}
    \phi_k^n(x^n,z^n) \approx \sum_{m=1}^{N_{ep}} \hat{\phi}^n_{k,m} \psi_{m}(x^n,z^n)=\sum_{i=1}^{N_{ep}} \phi_k^n(x_i^n,z_i^n) h_{i}(x^n,z^n), \quad \text{for} \quad (x^n,z^n) \in E^n,
\end{equation}
where the former and latter sum expressions are the modal and nodal forms, respectively. $N_{ep}$ is the number of local expansion coefficients, $\hat{\phi}^n_{k,m}$, or local solution values, $\phi_k^n(x_i^n,z_i^n)$. Moreover, $\psi_{m}$ is the modal basis functions (Jacobi polynomials) up to order $P$ and $h_{i}$ is the nodal basis functions (Lagrange polynomials) up to order $P$. The order of the basis functions, $P$, is connected to $N_{ep}$ as $N_{ep} =\frac{(P+1)(P+2)}{2}$. Each element is mapped onto a 2D simplex reference domain, $\mathcal{R}$, given by: $(r, s) \geq-1 ; r+s \leq 0$. Affine or curvilinear mappings, $(x,z) = \Psi(r,s)$ and $(r,s) = \Psi^{-1}(x,z)$, are introduced, ultimately giving rise to a local transformation Jacobian, $\mathcal{J}^n = \frac{\partial x^n}{\partial r} \frac{\partial z^n}{\partial s} - \frac{\partial x^n}{\partial s} \frac{\partial y^n}{\partial r}$. On $\mathcal{R}$, an orthonormal basis, $\psi_m(r,s)$, is constructed on a nodal set of grid points, $(r,s)$, as outlined in Hesthaven and Warburton (2007) \cite{hesthaven2007nodal}, using the family of Jacobi polynomials and Gauss-Lobatto points. Hereby leading to a well-behaved generalized Vandermonde matrix for interpolation
\begin{equation}
    \mathcal{V}_{ij} = \psi_j(r_i,s_i), \quad \text{for} \quad (i,j) = 1,...,N_{ep}.
\end{equation}

\subsubsection{Weak formulation and Galerkin spectral element discretization}

An essential step in the quest for the spatial discretization of the pseudo-impulsive formulation using the SEM is the derivation of the weak formulation in integral form. Ultimately, this results in the removal of second-order derivatives and enables the natural incorporation of Neumann boundary conditions into the scheme. Dividing the boundaries into Dirichlet and Neumann boundaries as $\Gamma = \Gamma^{DBC} + \Gamma^{NBC}$, and considering a test function on $\Omega$ as $v = v(x,z): \Omega \xrightarrow{} {\rm I\!R}$, that is assumed to be infinitely smooth and vanishing on $\Gamma^{DBC}$. By: i) multiplying the time-invariant Laplace problem in \eqref{eq:Laplace_equation} by $v$, ii) applying integration by parts, iii) using the properties of $v=0$ on $\Gamma^{DBC}$, and iv) rearranging, the weak formulation yields
\begin{equation}\label{eq:weak_formualtion}
         \int \int \nabla \phi_k \cdot \nabla v ~ dx dz = \int_{\Gamma^{NBC}} q ~ v ~ d \Gamma,
\end{equation}
where $q = q(x,y): \Gamma^{NBC} \xrightarrow{} {\rm I\!R}$ is the Neumann flux which is only relevant for inhomogeneous conditions on $\Gamma^{body}$ and $\Gamma^{\infty}$ in the case of a Sommerfeld condition. A Galerkin approach is adopted by setting the test function equal to the basis function, $v = \psi$, in the weak formulation, and then inserting the polynomial representation from \eqref{eq:polynomial_representation}. With this, a global linear system of equations of $N_\text{DOF}$ Degrees Of Freedom (DOF) is obtained in the form
\begin{equation}\label{eq:linear_system_of_equations}
    \mathcal{A} \boldsymbol{\phi_k} = \boldsymbol{b}, \quad \text{where} \quad \mathcal{A} \in \mathbb{R}^{N_\text{DOF} \times N_\text{DOF} }, \quad \text{and} \quad (\boldsymbol{\phi_k},\boldsymbol{b}) \in \mathbb{R}^{N_\text{DOF}},
\end{equation}
here $\mathcal{A}$ is the sparse system matrix, $\boldsymbol{b}$ is the system vector containing information about the boundary conditions, and $N_\text{DOF}$ being the number of nodal points in the spatial discretization. Local contributions to the system matrix, $\mathcal{A}^n$, can be derived\cite{engsig2016stabilised} as
\begin{equation}
    \mathcal{A}^n = \mathcal{D}_x^T \mathcal{M}^n \mathcal{D}_x + \mathcal{D}_z^T \mathcal{M}^n \mathcal{D}_z,
\end{equation}
with $\mathcal{M}^n$ being the local mass integration matrix, $\mathcal{M}^n = \mathcal{J}^n \mathcal{M} = \mathcal{J}^n \left ( \mathcal{V} \mathcal{V}^T \right ) ^{-1}$, and $\mathcal{D}_x$ and $\mathcal{D}_z$ is the differentiation matrices in the $x$- and $z$-direction, respectively, defined through the chain rule as $\mathcal{D}_x = \frac{\partial r}{\partial x} \mathcal{D}_{r} + \frac{\partial s}{\partial x} \mathcal{D}_{s}$ and $\mathcal{D}_z  = \frac{\partial r}{\partial z} \mathcal{D}_{r} + \frac{\partial s}{\partial z} \mathcal{D}_{s}$. Now, the differentiation operators on $\mathcal{R}$ are given as $\mathcal{D}_{r} = \mathcal{V}_{r} \mathcal{V}^{-1}$ and $\mathcal{D}_{s} = \mathcal{V}_{s} \mathcal{V}^{-1}$ with $\mathcal{V}_{r,ij} = \frac{\partial \psi_j}{\partial r} \mid _{(r_i,s_i)}$ and $\mathcal{V}_{s,ij} = \frac{\partial \psi_j}{\partial s} \mid _{(r_i,s_i)}$ for $(i,j) = 1,...,N_{ep}$. Similar contributions - this time in a 1D setting due to the line integral - can be derived for the system vector as
\begin{equation}
    \boldsymbol{b}^n = \boldsymbol{b}^n + \mathcal{M}^n_{1D} \boldsymbol{q}, \quad \text{with} \quad \boldsymbol{q} = [q_1,...,q_{N_{ep}^{1D}}]^T,
\end{equation}
where $\boldsymbol{b}$ is initialized to be all zeros, $\mathcal{M}^n_{1D}$ is the $n$'th mass integration matrix on the boundary element, and $\boldsymbol{q}$ are the discrete nodal flux values.

Inhomogeneous Dirichlet boundary conditions on $\Gamma^{FS}$ are not enforced due to the vanishing properties of the test function, $v$, on $\Gamma^{DBC}$. These are strictly imposed by modification of the system matrix and vector, $\mathcal{A}$ and $\boldsymbol{b}$, in standard ways.

\subsubsection{Curvilinear elements and non-constant transformation Jacobians}

For curved physical boundaries a geometrical second order error is introduced when using affine elements; ultimately, ruining the favorable spatial spectral convergence rate of the SEM. This challenge can be solved by considering curvilinear elements that will be constructed by modification of existing affine elements through transfinite interpolation with linear blending due to Gordon and Hall (1973) \cite{gordon1973transfinite,gordon1973construction}.

The curvilinear elements will - in contrast to the affine elements - result in non-constant transformation Jacobians leading to the possibility of aliasing errors in the approximation due to the nonlinear transformation in the discrete inner products, and as a result of this affecting the convergence rate. To overcome this challenge, a super-collocation method \cite{kirby2003aliasing,engsig2016stabilised,engsig2019mixed} is incorporated to handle the discrete inner products.


\subsection{Temporal integration}\label{sec:Time_Integration}

The method of lines approach leaves us with a semi-discrete system of ordinary differential equations from \eqref{eq:free_surface_conditions}, that are evolved in time using the classical ERK4. In combination with the time integration scheme, a global Courant Friedrichs Levy (CFL) condition is exploited to ensure sufficient stability, accuracy, and efficiency. Hereby, an upper bound for the time step size, $\Delta t$, is set to be
\begin{equation}\label{eq:time_step_size}
    \Delta t=\frac{\operatorname{Cr} \Delta x_{\min}}{u_{\max}},
\end{equation}
where $\operatorname{Cr} \in [0.5;1]$ is the Courant number, $\Delta x_{\min}$ is the minimum grid spacing on the free surface, and $u_{\max} = \sqrt{gh}$ is set to be the asymptotic shallow water limit. This choice of time step size can be seen to be quite conservative, as the propagation of higher frequency waves will need multiple time steps to move an equivalent $\Delta x_{\min}$. Using this time step size, the time domain, $\mathcal{T}$, is partitioned uniformly by $N+1$ sample points
\begin{equation}\label{eq:t_i}
    t_i = i \Delta t, \quad \text{for} \quad i = 0,1,...N. 
\end{equation}

\subsection{Computation of the pressure and hydrodynamic coefficients}

As is evident from \eqref{eq:F_jk}, the computation of the radiation force, $F_{jk}$, require information about the dynamic pressure, $p_k$. The evaluation of this quantity involves time derivatives of the radiation potentials that will be computed discretely by finite difference approximations \cite{leveque2007finite} as
\begin{equation}
    \frac{\partial \phi_k(x,z,t_i)}{\partial t} \approx \sum_{n = \alpha }^{\beta} a_n \phi(x,z,t_i+n\Delta t), \quad \text{for} \quad i = 0,1,...N,
\end{equation}
where the time vector, $t_i$, is discretized uniformly as in \eqref{eq:t_i} given the time step size from \eqref{eq:time_step_size}. The $a_n$ coefficients are determined though the \textit{method of undetermined coefficients} dependent of the order of accuracy and whether \textit{forward}, \textit{centered}, or \textit{backward} stencils are considered. To complement the $4$th order accurate time integration scheme, an equivalent accurate finite difference scheme is constructed using different stencils, and the convergence of the derivative approximation has been successfully validated.

The computation of the Fourier transforms needed in \eqref{eq:a_jk_b_jk} for the evaluation of the hydrodynamic coefficients is applied discretely through FFTs \cite{cooley1965algorithm}. The frequency domain, $\mathcal{H}$, is discretized uniformly with the frequency resolution, $\Delta f = 1 / t_{end}$, with $t_{end}$ being the time length of the simulation. To increase the resolution, the displacement and force signals, $x_k$ and $F_{jk}$, are extended with zeros assuming that the latter have reached a zero steady state solution, which will be ensured by following the design proposal in the upcoming section. This particular technique has been used extensively throughout the presented work.

\subsection{Construction of the discrete pseudo-impulsive displacement}\label{sec:discrete_pseudo_impulse}

On the body boundary, $\Gamma^{body}$, a pseudo-impulsive velocity, $\dot{x}_k$, is to be imposed. This velocity will be based on a Gaussian displacement signal, $x_k$, \cite{amini2017solving,amini2018pseudo} and designed such that: 1) the Gaussian is tailored to contain a specific range of frequencies, 2) the signal is practically zero at $t=0$, and 3) the resultant radiation force signal, $F_{jk}$, should reach a zero steady state solution at $t=T$. With these requirements, a unit height Gaussian displacement signal is defined as
\begin{equation}\label{eq:pseudo_impulse}
    x_{k}(t)=e^{-2 \pi^{2} s^{2}\left(t-t_{0}\right)^{2}}, \quad \text{with} \quad  t_0 = \sqrt{\frac{\log(\epsilon)}{-2 \pi^2 s^2}}, \quad \text{and}, \quad s=\sqrt{\frac{-f_r^{2}}{2 \ln (r)}},
\end{equation}
where $t_0$ is the time location of the Gaussian peak and $s$ is a frequency range governing parameter. The parameter, $\epsilon$, is some small quantity to ensure that $x_k(0) \approx 0$ and $r$ is governing how many orders of magnitude there should be between the frequency components $\hat{x}_k(f_r)$ and $\hat{x}_k(0)$, where $f_r$ is the chosen maximum resolvable design frequency so that the Gaussian contains frequencies between $0$ and $f_r$. For this work, $r=10^{-4}$ has been used. The maximum resolvable design resolution is defined through linear wave theory as
\begin{equation}\label{eq:design_wave}
    f_r = \frac{\omega_r}{2 \pi}, \quad \omega_r=\sqrt{g k_{r} \tanh \left( k_{r} h \right )}, \quad k_r = \frac{2 \pi}{L_r}, \quad \text{and} \quad L_r = \alpha ~ \Delta x_{\max},
\end{equation}
with $\Delta x_{\max}$ being the largest distance between two nodal points on the free surface boundary, $\Gamma^{FS}$, and $\alpha$ is a parameter governing the length of the smallest resolvable design wavelength. Unless stated otherwise, $\alpha = 3$, resulting in $4$ grid points per smallest possible wavelength. Analytical expressions for the pseudo-impulsive velocity and the Fourier transform of the displacement signal can be obtained as:
\begin{equation}
     \dot{x}_k(t) = -4 \pi^2 s^2 (t-t_0) e^{-2\pi^2s^2(t-t_0)^2}, \quad \text{and} \quad \mathcal{F}\{x_k(t)\} = \hat{x}_k(f) = \frac{1}{s \sqrt{2 \pi}} e ^{\frac{-f^2}{2s^2}}.
\end{equation}
Finally, the simulation time is set to be at least $t_{end} = 3 t_0$ to ensure that $F_{jk}(t_{end}) \approx 0$.
\section{Numerical Properties}\label{sec:Numerical_Properties}

In the following, the numerical properties of the developed model are highlighted. This will be done by considering: 1) the spatial accuracy of the Laplace solver in terms of a convergence study, 2) the efficiency when solving the sparse linear system of equations in \eqref{eq:linear_system_of_equations}, and 3) the temporal stability in relation to the spatial discretization.

\subsection{Convergence of Laplace solver}

As an exact solution to the pseudo-impulsive radiation problem is non-trivial, we consider the convergence study of the Laplace solver in the scope of the Method of Manufactured Solutions (MMS) \cite{roache2002code}, where a true solution, $\phi_k^{MMS}$, is assumed to the problem, from which analytical expressions can be derived for the boundary conditions and right hand side function, and hereby reforming the Laplace problem into a Poisson problem. The true solution should be infinitely smooth, e.g. a combination of trigonometric functions, to complement the wave problem, and with this, we seek to verify that we are able to archive: 1) spectral $P$-convergence for fixed meshes by increasing the order of orthonormal basis functions, $\psi$, and 2) algebraic $h$-convergence of order $P+1$ for fixed orders of orthonormal basis functions with decreasing element size. 

\begin{figure}[t]
\begin{minipage}[c]{0.3\textwidth}
    \caption{Three unstructured uniform meshes using triangular elements with decreasing element size, $N_{elm} = (32,53,73)$, for discretization of a floating cylinder with radius $R$.}
    \label{fig:First_three_meshes_affine}
\end{minipage}
\hfill
\begin{minipage}[c]{0.69\textwidth}
    \includegraphics{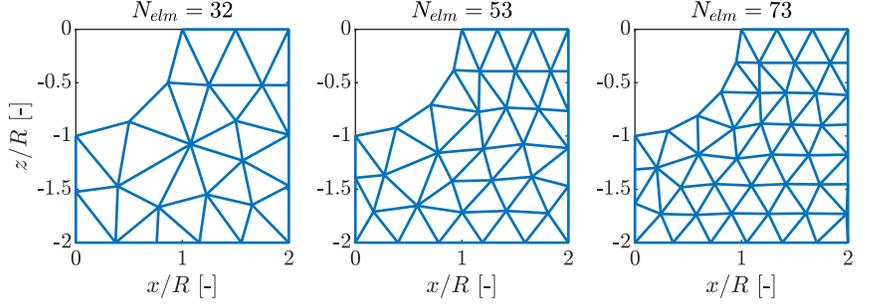}
\end{minipage}
\end{figure}

\begin{figure}[t]
\begin{minipage}[c]{0.3\textwidth}
    \caption{Nodal distribution (red dots) of affine and curvilinear elements near $\Gamma^{body}$ for different polynomial orders, $P=(1,2,3)$, on the coarsest mesh in Figure \ref{fig:First_three_meshes_affine}, $N_{elm}=32$.}
    \label{fig:Mesh_curvilinear}
\end{minipage}
\hfill
\begin{minipage}[c]{0.69\textwidth}
    \includegraphics{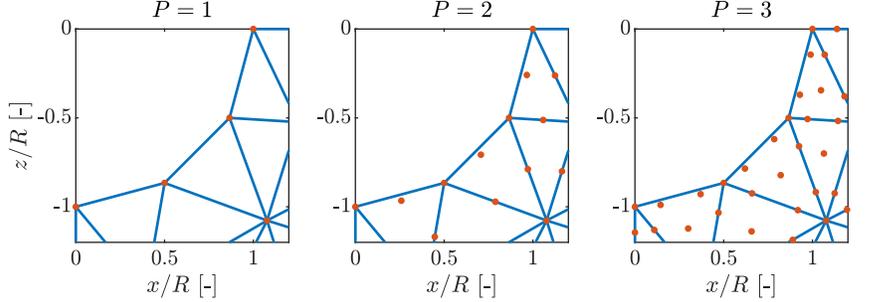}
\end{minipage}
\end{figure}

For the convergence study, we consider a fluid domain, $\Omega$, that - due to the symmetry boundary - has a quarter of a floating cylinder with radius, $R$, in the top left corner. The domain is partitioned using curvilinear elements around $\Gamma^{body}$ to capture the curved boundary of the structure, and for the remaining part, we seek to use regular affine elements in an unstructured fashion. When considering the convergence study, the elements are arranged uniformly, as seen in Figure \ref{fig:First_three_meshes_affine} for the three coarsest meshes used having a total of $N_{elm} = (32,53,73)$ number of elements, respectively. In Figure \ref{fig:Mesh_curvilinear}, a zoomed-in plot of the body boundary elements for the mesh with $N_{elm} = 32$ is pictured to show how the distribution of nodal points are placed within the curvilinear elements for $P = (1,2,3)$.

\begin{figure}[t]
\begin{minipage}[c]{0.3\textwidth}
    \caption{$P$-convergence study using curvilinear and affine unstructured triangular elements. Fixed meshes from Figure \ref{fig:First_three_meshes_affine} with increasing polynomial order, $P$.}
    \label{fig:P_Convergence_curvilinear}
\end{minipage}
\hfill
\begin{minipage}[c]{0.60\textwidth}
    \includegraphics{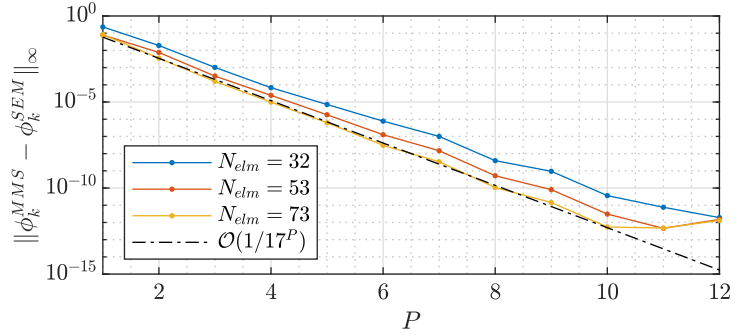}
\end{minipage}
\end{figure}

\begin{figure}[t]
\begin{minipage}[c]{0.3\textwidth}
    \caption{$h$-convergence study using curvilinear and affine unstructured triangular elements. Fixed polynomial order, $P$, and decreasing elements size in terms of maximal element edge length, $h_{max}$.}
    \label{fig:h_Convergence_curvilinear}
\end{minipage}
\hfill
\begin{minipage}[c]{0.60\textwidth}
    \includegraphics{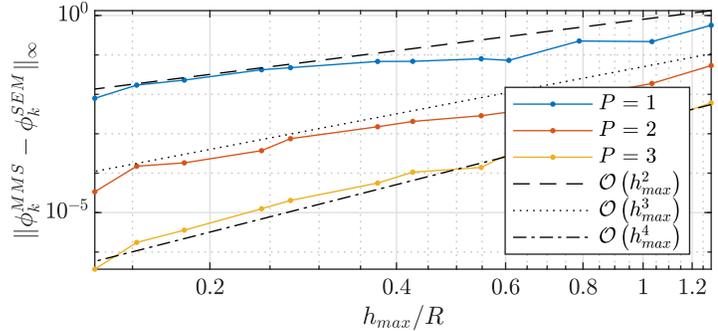}
\end{minipage}
\end{figure}

The convergence study is performed as shown in Figure \ref{fig:P_Convergence_curvilinear} and Figure \ref{fig:h_Convergence_curvilinear} with the $P$- and $h$-test, respectively, where the error is evaluated globally in the solution to the partial differential equation using the infinity norm as $\|\phi_k^{MMS} - \phi_k^{SEM} \|_{\infty}$. For the former test, the three meshes in Figure \ref{fig:First_three_meshes_affine} are used with increasing order of the orthonormal basis functions, $P=1,...,12$. The spectral convergence behavior is confirmed in the semi-logarithmic plot by comparing with a $\mathcal{O}(1/17^P)$ line, where the error can be seen to reach about $\approx 10^{-13}$, from where round-off errors start to become significant and accumulate. It can also be noted that the finer meshes, in general, give smaller errors than the coarser ones. This is an early indication of $h$-convergence, which ultimately can be confirmed in the latter test on the logarithmic plot, showing algebraic convergence of order $P+1$ for $P=(1,2,3)$ by decreasing the size of the elements.

\subsection{Scaling of the computational effort}

When solving free surface potential flow, the main concern and computational bottleneck is the work associated with solving the Laplace problem numerically, which has to be done multiple times - dependent on the time discretization method - at every time step. Therefore, one should seek to complete this challenge as efficiently as possible.

\begin{figure}[b]
\begin{minipage}[c]{0.3\textwidth}
    \caption{Analysis of computational work effort for the 2D SEM solver using a sparse reverse Cuthill-McKee ordering. Each simulation (black dot) represent an average of solving the Laplace problem $10$ times.}
    \label{fig:CPU}
\end{minipage}
\hfill
\begin{minipage}[c]{0.60\textwidth}
    \includegraphics{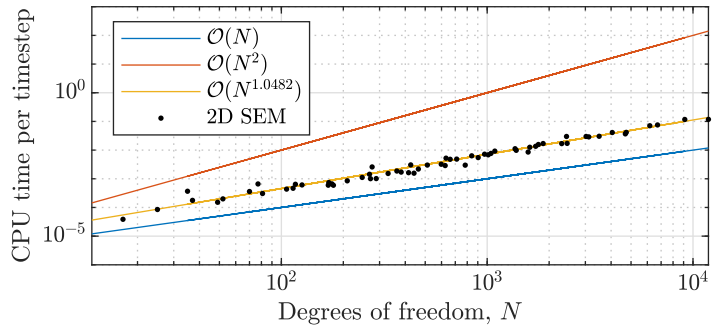}
\end{minipage}
\end{figure}

The linear system of equation in \eqref{eq:linear_system_of_equations} is solved directly in MATLAB, $\verb|A\b|$, using a sparse reverse Cuthill-McKee reordering \cite{CuthillMckee}. The permutation, $q$, reorders the system matrix as $\mathcal{A}(q,q)$, such that the band width is minimized, hence making the direct solver much more efficient. In Figure \ref{fig:CPU}, the study of computational effort can be seen in terms of the DOF of the system, $N_\text{DOF}$, against CPU time per time step. Nine different uniform meshes where constructed with the ratio between the coarsest and the finest mesh as $h_{coarse} / h_{fine} \approx 7$ yielding meshes having $N_{elm} \in [17;365]$ number of elements. For each mesh, polynomial orders of $P=1,...,8$ were applied, resulting in $N_\text{DOF} \in [17;11885]$. From the study, $\mathcal{O}(N^p)$ can be observed with $p = 1.0482$, when solving the linear system using a direct solver. The reason why $p \neq 1$ is due to fill-in of $\mathcal{A}$ by the direct solver.

Despite the fact that the sparse scheme - due to the local support of the basis functions - already is very efficient, the computational effort can be made to scale with $\mathcal{O}(N)$ by considering the use of a geometric $p$-multigrid method that exploits the polynomial convergence features of the SEM model to solve the Laplace problem iteratively \cite{laskowski2019modelling,engsig2021efficient}.

\subsection{Linear stability analysis of spatial discretization}

Temporal instability is not solely related to temporal discretization as outlined in Robertson and Sherwin (1999) \cite{robertson1999free} where the issues related to computational instabilities of incompressible and inviscid free surface flow are discussed in the context of the SEM. The reason was found to be caused by asymmetrical spatial discretizations of the free surface of the fluid domain, for which the authors proposed a method for studying the stability of the scheme. Due to the adopted MOL approach, we are left with a semi-discrete system for the kinematic and dynamic free surface condition from \eqref{eq:free_surface_conditions}. The eigenvalues and eigenvectors of this system will dictate the temporal stability of the scheme, where eigenvalues with positive real components indicate an unstable system, as these govern the diffusive properties, thus causing the system to diverge. Such analysis is performed in a recent work \cite{engsig2019mixed}.

One way to ensure stability when solving free surface flows using the SEM is by applying quadrilateral elements for the entirety of the domain \cite{engsig2016stabilised,mortensen2021simulation}; however, in Engsig-Karup et al. (2019) \cite{engsig2019mixed} it was shown to be sufficient to exploit a single layer of this type of elements below the free surface boundary, $\Gamma^{FS}$. It was found that the instability issue is related to the accuracy of computing the vertical free surface velocity, which govern the dispersive property of the discrete model. To maximize this, the nodal distribution beneath $\Gamma^{FS}$ must be on a vertical line. In the same paper, weak stability properties, when using triangular elements, were highlighted, showing that some polynomial orders were keen to provide instability issues despite the fulfilled symmetry requirement. Similar analyses have been carried out in relation to this work, and similar results have been obtained. Hereby emphasizing the fact that stability analyses are solely related to the specific meshes; with this, one should be careful - when using triangular elements - picking the proper polynomial order that allows for a stable temporal evolution of the governing free surface equations.

%

\section{Numerical Experiments}\label{sec:Numerical_Experiments}

We seek to outline simulation results for floating bodies (a cylinder and a box) in the following by considering comparisons of force signals and added mass and damping coefficients with other well-established numerical models and benchmark results from the literature. Also, an extensive investigation of the numerical phenomenon of free surface \textit{spurious oscillations} arising from unresolved energy in the pseudo-impulse is carried out.

\subsection{Heaving and surging structures}

\begin{figure}[t]
\begin{minipage}[c]{0.3\textwidth}
    \caption{Visualization of mesh used for simulation of a floating cylinder. Mesh refinement on $\Gamma^{FS}$ and $\Gamma^{body}$ with 5 curvilinear elements constituting $\Gamma^{body}$ for meshing of a quarter of the cylinder due to symmetry in the problem. The radius is $R = 1$ m, the domain length is $L = 10$ m, and the finite water depth is $h \approx 6.283 $ m.}
    \label{fig:Mesh_Cylinder}
\end{minipage}
\hfill
\begin{minipage}[c]{0.69\textwidth}
    \includegraphics{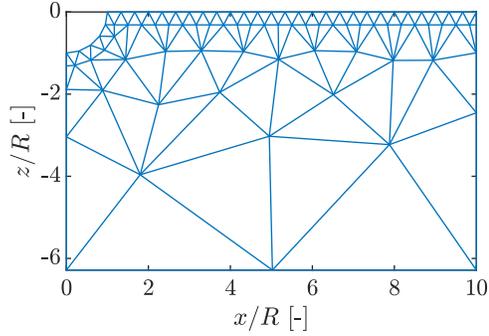}
\end{minipage}
\end{figure}

For the simulations presented in the following, affine and curvilinear elements have been used, where the elements on the free surface boundary, $\Gamma^{FS}$, is constructed symmetrically to stabilize most of the polynomial orders. Only stable orders have been considered. The elements will be refined toward the free surface and body boundaries, $\Gamma^{FS}$ and $\Gamma^{body}$, e.g. as shown in Figure \ref{fig:Mesh_Cylinder}. On the far-field boundary, $\Gamma^{\infty}$, a combination of a Sommerfeld condition \cite{zhang2009non} and an absorption zone technique \cite{larsen1983open} are used to make the domain non-reflective for this multi-frequency wave problem.

\subsubsection{The floating cylinder}
\begin{figure}[t]
\begin{minipage}[c]{0.3\textwidth}
    \caption{Comparisons of heave-heave radiation force signals, $F_{33}(t)$, of the FDM-based and the SEM models. Simulation obtained using $P=3$ ordered basis functions on the mesh shown in Figure \ref{fig:Mesh_Cylinder}. The design wave period for the displacement signal is denoted by $T$.}
    \label{fig:Comparison_Force}
\end{minipage}
\hfill
\begin{minipage}[c]{0.69\textwidth}
    \includegraphics{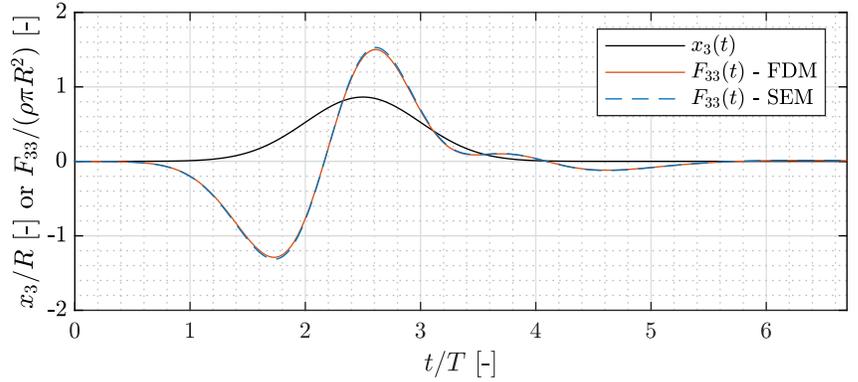}
\end{minipage}
\end{figure}

The first floating structure is the classical half submerged cylinder with a radius $R$. We compare the calculations with an equivalent FDM-based solver, used to obtain the results in Read and Bingham (2012) \cite{read2012overset} and others \cite{read2011linear,read2012solving}, by exposing the two models to identical input displacement signals, $x_k(t)$, and then studying the corresponding output force signals, $F_{jk}(t)$. For the simulation, we used the mesh shown in Figure \ref{fig:Mesh_Cylinder}, having 5 curvilinear elements around $\Gamma^{body}$ for meshing of a quarter of the cylinder surface. Basis functions of order $P=3$ are chosen to ensure stability. From the comparison in Figure \ref{fig:Comparison_Force}, an excellent visual agreement can be observed as the two force signals are placed on top of each other yielding non-dimensional absolute errors in the order of $\mathcal{O}(0.01)$. The minor discrepancies are located around the peaks of the force signals and are associated with differences in spatial resolution.

Next, a comparison of added mass and damping coefficients with results published in Read and Bingham (2012) \cite{read2012overset}. These benchmark results have been validated against analytical infinite-depth solution calculated using a multipole method \cite{ursell1949heaving}. The radius of the cylinder is set to $R = 0.5$ m and the water depth is $h = 3$ m yielding $R/h = 6$ [-]. The added mass and damping coefficients have been normalized as
\begin{equation}
    \mu_{jk} = \frac{a_{jk}}{\frac{1}{2} \pi \rho R^2 }, \quad \nu_{jk} = \frac{b_{jk}}{\frac{1}{2} \pi \rho \omega R^2 }, \quad \text{for} \quad (j,k) = 1,3.
\end{equation}

\begin{figure}[b]
    \begin{minipage}[t]{.45\textwidth}
        \centering
        \includegraphics[width=\textwidth]{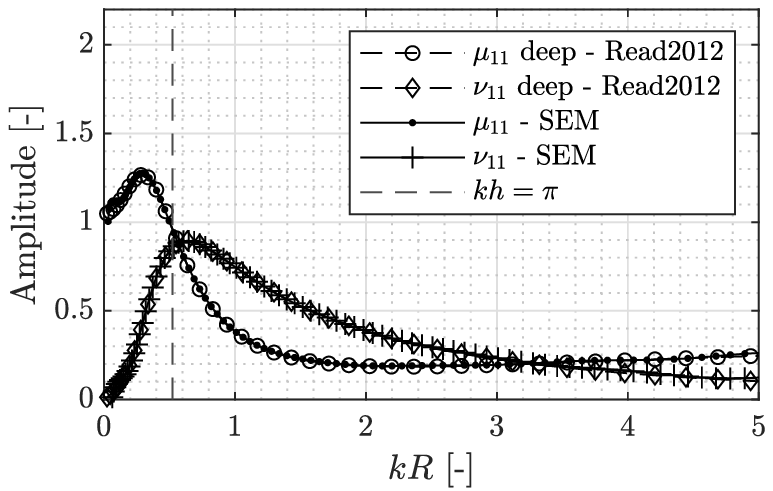}
        \subcaption{Surge-Surge case: $\mu_{11}$ and $\nu_{11}$.}
    \end{minipage}
    \hfill
    \begin{minipage}[t]{.45\textwidth}
        \centering
        \includegraphics[width=\textwidth]{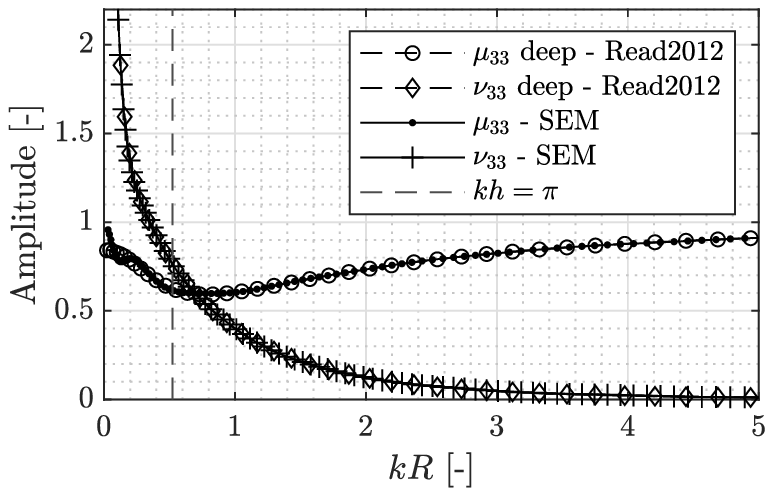}
        \subcaption{Heave-Heave case: $\mu_{33}$ and $\nu_{33}$.}
    \end{minipage}  
    \caption{Dimensionless added mass and damping coefficients for a floating cylinder of radius $R = 0.5$ m and a finite water depth of $h = 3$ m. Basis functions of order $P=3$ is used. Comparison between FDM-based solver results \cite{read2012overset,read2012solving} and SEM model.}
    \label{fig:Cylinder_P3_Read2012}
\end{figure}

The surge-surge and heave-heave results are displayed in Figure \ref{fig:Cylinder_P3_Read2012} (left) and  \ref{fig:Cylinder_P3_Read2012} (right), respectively. The numerical results show close to perfect visual agreement, especially for $kh > \pi$ where the waves - in terms of wave length - are sufficiently short to be uninfluenced by the finite depth. The latter comparison indeed highlights the legitimacy of the proposed model.

\subsubsection{The floating box}
The second floating structure to consider is a floating box-type structure with the draft, $d$, and the half-length, $a$. A comparison with the analytical method in an infinite fluid domain of finite water depth presented by Zheng et al. (2004) \cite{zheng2004radiation} is to be carried out. The aforementioned method was validated against Lee (2005) \cite{lee1995heave} and a BEM model. 

For the simulation, we adopt the same setting as in the comparison study, hereby using a dimensionless depth of $h/d = 3$ [-] and dimensionless length of $a/d = 0.5$ [-]. As all boundaries are non-curved, solely affine elements will be considered. Regarding the meshing of the box, 6 elements constitutes the bottom of the structure and 12 elements the side. Basis functions of order $P = 4$ is used. The added mass and damping coefficients have been normalized as
\begin{equation}
    \mu_{jk} = \frac{a_{jk}}{2 \rho a d }, \quad \nu_{jk} = \frac{b_{jk}}{2  \rho \omega a d }, \quad \text{for} \quad (j,k) = 1,3.
\end{equation}

\begin{figure}[t]
    \begin{minipage}[t]{.45\textwidth}
        \centering
        \includegraphics[width=\textwidth]{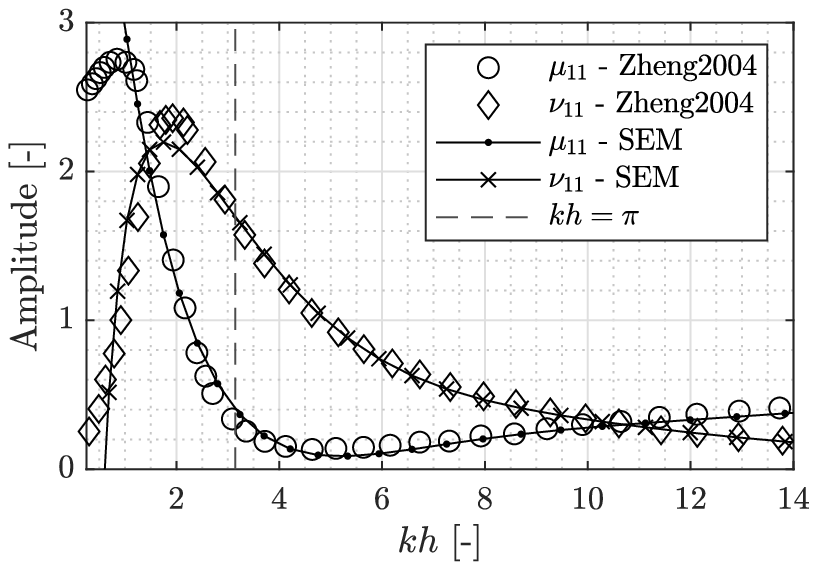}
        \subcaption{Surge-Surge case: $\mu_{11}$ and $\nu_{11}$.}
    \end{minipage}
    \hfill
    \begin{minipage}[t]{.45\textwidth}
        \centering
        \includegraphics[width=\textwidth]{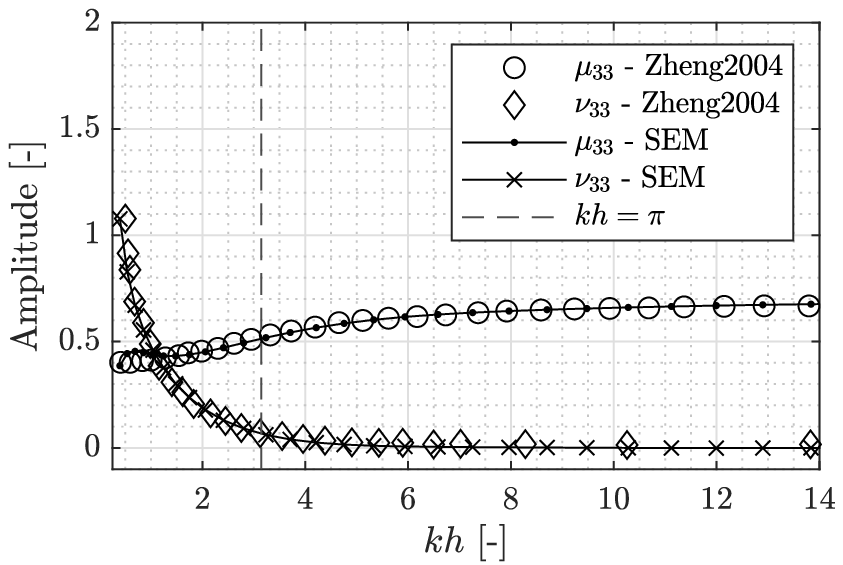}
        \subcaption{Heave-Heave case: $\mu_{33}$ and $\nu_{33}$.}
    \end{minipage}  
    \caption{Dimensionless added mass and damping coefficients for a floating box with dimensionless depth, $h/d = 3$ [-], and dimensionless length, $a/d = 0.5$ [-]. Basis functions of order $P=4$ is used. Comparison between Zheng et al. (2004) \cite{zheng2004radiation} and the SEM model.}
    \label{fig:Box_P4_Zheng2004}
\end{figure}

In Figure \ref{fig:Box_P4_Zheng2004} (left) and (right) the numerical results for non-dimensional added mass and damping coefficients can be seen for surge-surge and heave-heave, respectively. For $kh > \pi$, an excellent visual agreement can be discovered as for the cylinder case. Note that the results from Zheng et al. (2004) have been manually digitized. 

\subsection{Investigation of spurious oscillations}\label{sec:Investigation}

As presented in the introduction, the difference between the classical impulsive formulation \cite{korsmeyer1988first,king1988seakeeping,bingham1994simulating} and the pseudo-impulsive formulation \cite{read2011linear,read2012overset,read2012solving,amini2017solving,amini2018pseudo} will be highlighted in the following. Starting with the former, we have that the impulsive motion is imposed rapidly at initial time, $t=0$, using the generalized Dirac Delta function, $\dot{x}_k(t) = \delta(t)$, i.e. a pressure release. From the Fourier transformation, we conclude that $\mathcal{F}\{ \delta(t) \} =  \hat{\delta}(f) = 1$, hence unit amplitude energy is imposed on all frequencies in the spectrum. This operation works well in the continuous setting; however, when the governing formulation is discretized, one faces the challenge of resolving energy on frequency components that are unresolvable by the discrete scheme. From experiments with the classical impulsive formulation, the unresolved energy has been seen be alias towards lower frequencies and give rise to \textit{spurious oscillations} of the free surface quantities. 

This is the primary motivation for using the pseudo-impulsive formulation, where the frequency spectrum is tailored through a unit-amplitude Gaussian displacement function as outlined in Section \ref{sec:discrete_pseudo_impulse}. From the experiences with the classical formulation, we seek to study the unresolved energy in terms of the spurious free surface oscillations using the pseudo-impulsive formulation as this - to the authors' knowledge - has not been investigated before. The study will be carried out by including too high frequencies; that will enable us to explore the behavior of the system. For the analysis, we solely consider the case of a floating cylinder in a heave-heave motion and compare it with a BEM-based model on infinite depth. The BEM-model has been validated against solutions on infinite-depth computed using the method of Porter (2008) \cite{porter2008solution}. Also, a stable polynomial order of $P=4$ is used on a mesh with - until further notice - $5$ curvilinear elements around a quarter of the cylinder surface. The finite water depth is set to $h/R = 6$ [-].

\subsubsection{Pseudo-impulse response}

\begin{figure}[t]
    \centering
    \includegraphics{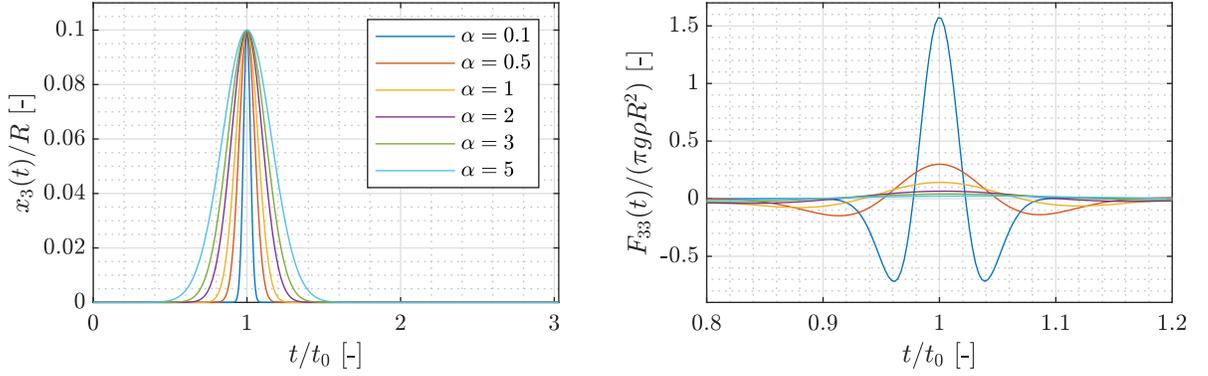}
    \caption{Left: Dimensionless displacement signals, $x_3(t)$, for $\alpha \in [0.1; 5]$. Right: Dimensionless force signals, $F_{33}(t)$, for $\alpha \in [0.1; 5]$.}
    \label{fig:SO_Displacement_and_Force}
\end{figure}

We first investigate the response when the same discrete system is exposed to different pseudo-impulses. Recall the parameter, $\alpha$, from \eqref{eq:design_wave} governing the smallest resolvable wavelength, $L_r = \alpha \Delta x_{\max}$, where $\Delta x_{\max}$ is the largest distance between two grid points on the free surface. For the results presented so far, a value of $\alpha = 3$ has been used with the argument that a linear wave requires $\pi$ (approximately $4$) points to be represented adequately. Thus, we arrive at $L_r = 3 \Delta x_{\max}$. In the following, we change the $\alpha$-parameter in the range $\alpha \in [0.1; 5]$, yielding both under- and over-resolved energy spectra. All of the displacement signals, $x_3(t)$, can be seen on Figure \ref{fig:SO_Displacement_and_Force} (left), where the same peak location, $t_0 = 8$ s, have been used. The corresponding force signals, $F_{33}(t)$, can be seen on Figure \ref{fig:SO_Displacement_and_Force} (right). Observe how decreasing $\alpha$ implies larger gradients in the displacement signal, hence larger boundary imposed velocities, ultimately yielding higher forces. Recall that a narrow displacement signal implies a wide band of frequencies, and vice versa for a wide signal. 

\begin{figure}[t]
    \centering
    \includegraphics{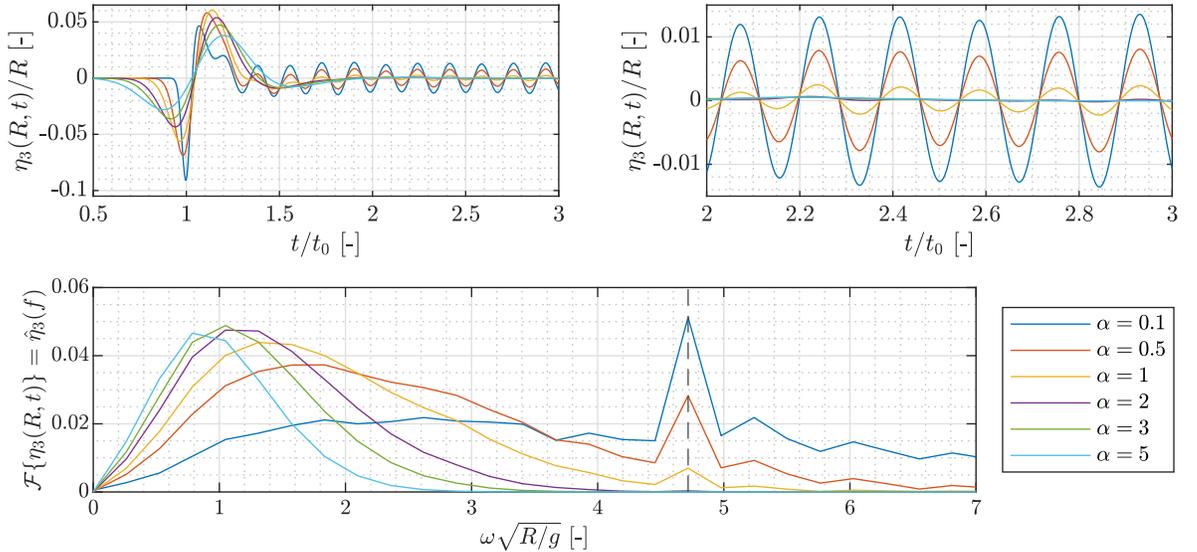}
    \caption{Top left: Dimensionless free surface elevation, $\eta_3(R,t)$, at the point between $\Gamma^{body}$ and $\Gamma^{FS}$ for $\alpha \in [0.1; 5]$. Top right: Zoomed-in snapshot of the top left plot. Bottom: Fourier transformation of $\eta_3(R,t)$ for $\alpha \in [0.1; 5]$.}
    \label{fig:SO_Eta_v1}
\end{figure}

When studying the free surface elevation, $\eta_3(R,t)$, right at the point between $\Gamma^{body}$ and $\Gamma^{FS}$, $(x,z) = (R,0)$, the consequence of using different pseudo-impulses can be seen in Figure \ref{fig:SO_Eta_v1} (top left and right). Here the presence of the spurious oscillations starts to take place for $\alpha < 2$, which confirms the initial use of $\alpha = 3$ to match the smallest resolvable wavelength. Also, the occurrence of the oscillations can be seen to be of constant phase independent of the pseudo-impulse; yet, the amplitudes tend to increase for decreasing $\alpha$. On Figure \ref{fig:SO_Eta_v1} (bottom), the Fourier transformation of the free surface signal, $\hat{\eta}_3(f) = \mathcal{F}\{ \eta_3(R,t) \}$, is shown, ultimately confirming the aforementioned observation of same phase and amplitude when $\alpha$ decreases as seen in increased spectral energy at the same particular frequency, $\omega \sqrt{R/g} \approx 4.7$. With this, it can be concluded that the spatial discretization governs which frequency the unresolved velocity imposed energy is aliased towards. Furthermore, it can again be verified that $\alpha \geq 3$ seems to provide the right physical results. \textbf{Remark:} Identical behavior have been observed in terms of the free surface velocity potential, $\phi_3(R,0,t)$.

\begin{figure}[t]
\begin{minipage}[b]{0.3\textwidth}
    \caption{Top: Dimensionless added mass coefficients, $a_{33}$, for $\alpha \in [0.1; 5]$. Bottom: Dimensionless damping coefficients, $b_{33}$, for $\alpha \in [0.1; 5]$. Both compared with an infinite-depth BEM model.}
    \label{fig:SO_Coefficients_v1}
\end{minipage}
\hfill
\begin{minipage}[c]{0.69\textwidth}
    \includegraphics{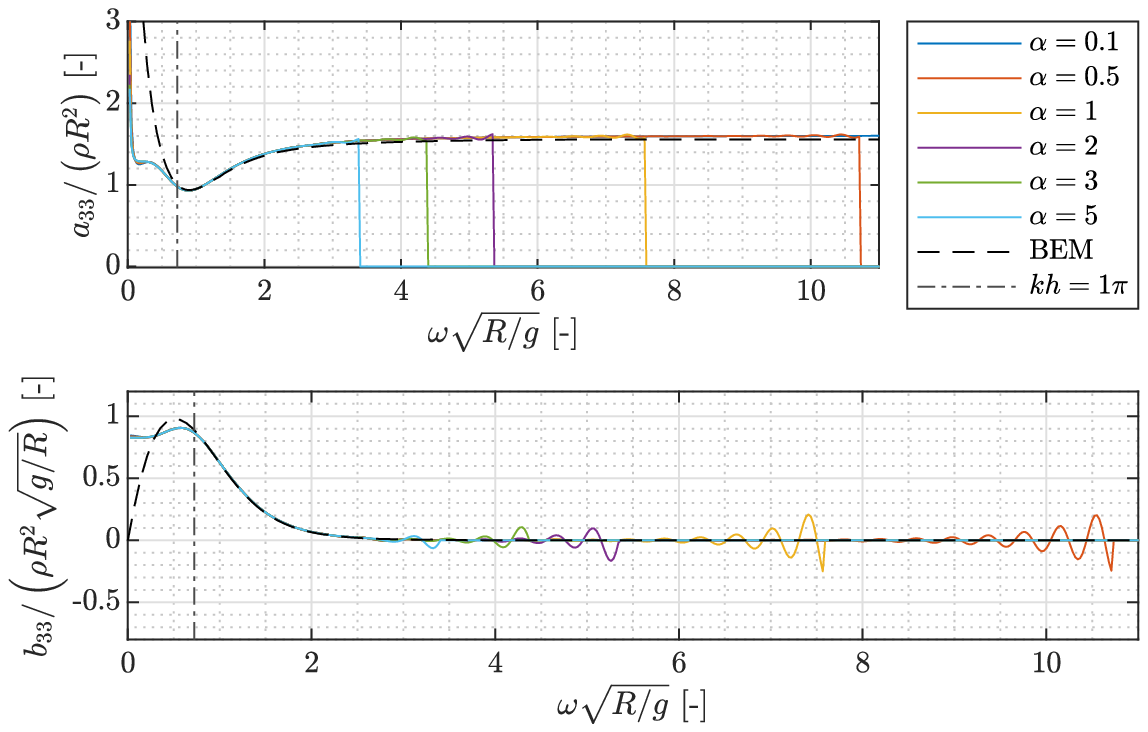}
\end{minipage}
\end{figure}

On Figure \ref{fig:SO_Coefficients_v1}, the dimensionless added mass and damping coefficients, $a_{33}$ and $b_{33}$, are plotted for $\alpha \in [0.1; 5]$. It can be observed how narrowing the displacement signal yields an increased frequency span. Note that the results for $\alpha = 0.1$ goes out to $\omega_r \sqrt{R/g} \approx 24$, hence not shown on the figure. In general, it is evident that near the cut-off frequency, $\omega_r\sqrt{R/g}$, some other oscillations are present. These are caused by the fact that the added mass and damping coefficient are computed via: $\mathcal{F} \{F_{jk}(t)\}/\mathcal{F}\{x_k(t)\}$. Hence near $\omega_r \sqrt{R/g}$, a small number is divided by another small number, thus prone to give rise to this kind of disturbance. Another result is how the added mass and damping coefficients can be well-approximated by decreasing the value of $\alpha$ - even to a level with the smallest introduced waves being far from adequately resolved by the spatial scheme. Somehow, the non-physical numerical time-domain spurious oscillations do not impact the frequency solution. The apparent advantage of this observation is that one can evaluate added mass and damping coefficients extremely efficiently by choosing sufficiently small values of $\alpha$, resulting in a narrow displacement signal, requiring less simulation time and a smaller domain size. \textbf{Remark:} The efficient computational approach was tried for the FDM-based pseudo-impulsive model \cite{read2011linear,read2012overset,read2012solving}, and similar result was discovered.

\subsubsection{Spatial resolution}

The next to investigate is the dependency on spatial resolution. Choosing the displacement signal to be very narrow by setting $s=1$ in \eqref{eq:pseudo_impulse} will enforce oscillations of the free surface quantities. Then, we seek to apply different meshes with increasingly more elements at $\Gamma^{FS}$ and around $\Gamma^{body}$. Let the parameter, $\beta$, govern the number of elements on a quarter of the cylinder, where the element edge length on $\Gamma^{FS}$ is (approximately) the same as on $\Gamma^{body}$. Recall from the previous investigation that $\beta = 5$ was used, whereas $\beta \in [3; 10]$ is considered in the following. 

\begin{figure}[t]
    \centering
    \includegraphics{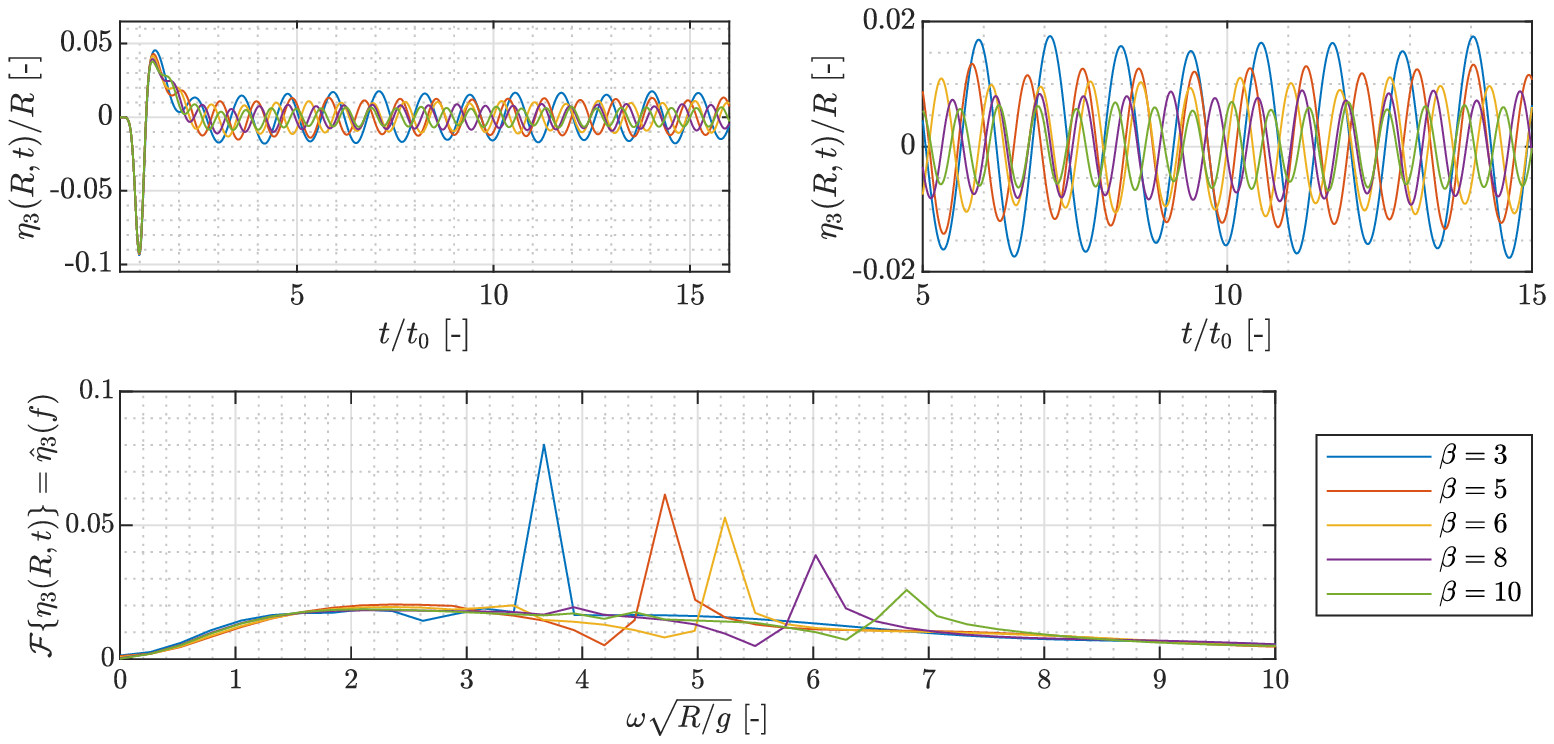}
    \caption{Top left: Visualization of free surface elevation, $\eta_3(R,t)$, at the point between $\Gamma^{body}$ and $\Gamma^{FS}$ for $\beta \in [3; 10]$. Top right: Zoomed-in snapshot of the top left plot. Bottom: Fourier transformation of $\eta_3(R,t)$ for $\beta \in [3; 10]$.}
    \label{fig:SO_Eta_v2}
\end{figure}

Due to the identical displacement signals, the force signals are - to a large extend - identical; yet, minor discrepancies will occur due to resolution. When considering the free surface elevation, $\eta_3(R,t)$, shown in Figure \ref{fig:SO_Eta_v2} (top left and right), the spurious oscillations are very much present for all simulations; however, the phases are not equal for different $\beta$, and the amplitude can be seen to decrease for increasing $\beta$. Also, the frequency seems to increase for increasing $\beta$. In Figure \ref{fig:SO_Eta_v2} (bottom) the equivalent frequency content of the free surface elevation, $\hat{\eta}_3(f)$, can be seen. Here, the time-domain observations can be confirmed, as it is clear that by increasing the spatial resolution: 1) the energy of the spurious oscillation is decreased, and 2) the location of the oscillation frequency component migrates towards higher frequencies. With this, it can be concluded that by refining the spatial resolution in the presence of unresolved energy, the frequency, on which the energy is aliasing towards, becomes higher and less significant.

\section{Conclusion}\label{sec:Conclusion}

A novel 2D spectral element model is proposed for the simulation of linear radiation potentials in relation to floating offshore structures. The mathematical formulation was established within the framework of linear potential flow using a pseudo-impulsive time domain formulation. The numerical methods and solution approaches were established for solving the governing equations numerically, where the spectral element method is handling the discretization of the spatial domain. 

The numerical properties of the model were tested in terms of spatial accuracy, computational efficiency, and temporal stability, with this confirming $P$- and $h$-convergence, $\mathcal{O}(N^{p})$ scalability, where $p \approx 1$, and instability considerations, respectively. At last, various numerical experiments were carried out and compared against established solvers and benchmark results from the literature. All showing excellent visual agreements and ultimately confirming the legitimacy of the proposed model. Also, a thorough novelty investigation of non-resolved velocity imposed energy leading to spurious oscillations of the free surface quantities was performed. With this, concluding and emphasizing the importance of correctly tailoring the Gaussian input signal, and also discovering an efficient non-physical way of computing added mass and damping coefficients. 

\bibliographystyle{unsrtnat}
\bibliography{references}  


\end{document}